\documentclass[preprint,12pt]{elsarticle}

\usepackage{amsmath}
\usepackage{epsfig}
\usepackage{amsmath}
\usepackage{amssymb}
\usepackage[mathscr]{euscript}
\usepackage{graphicx}
\usepackage{color}
\usepackage{algpseudocode}
\usepackage{algorithm}
\usepackage{pstricks-add}
\usepackage{pst-tree}
\usepackage{multirow}
\usepackage{rotating}

\newcommand{\newendproof}{\hspace{10pt}\rule{6pt}{6pt}}

\newcommand{\bee}{\begin{eqnarray}}
\newcommand{\eee}{\end{eqnarray}}
\def\uun{{\bf 1}}
\def\un{{\bf 0}}

\newcommand{\erre}{{\mathbb R}}

\newcommand{\eps}{\varepsilon}

\newcommand{\pf}{\noindent{\bf Proof\ \ }}

\def\eps{{\epsilon}}

\newtheorem{prop}{Proposition}[section]
\newtheorem{thm}{Theorem}[section]
\newtheorem{lemma}{Lemma}[section]
\newtheorem{cor}{Corollary}[section]

\newtheorem{defini}{Definition}[section]

\journal{Games and Economic Behavior}

\begin{document}

\begin{frontmatter}

\title{Strong Nash equilibria and mixed strategies}

\author[polimiMATE]{Eleonora Braggion}
\author[polimiDEIB]{Nicola Gatti}
\author[polimiMATE]{Roberto Lucchetti}
\author[CMU]{Tuomas Sandholm}
\address[polimiMATE]{Dipartimento di Matematica, Politecnico di Milano, piazza Leonardo da Vinci 32, 20133 Milano, Italy}
\address[polimiDEIB]{Dipartimento di Elettronica, Informazione e Bioningegneria, Politecnico di Milano, piazza Leonardo da Vinci 32, 20133 Milano, Italy}
\address[CMU]{Computer Science Department, Carnegie Mellon University, 5000 Forbes Avenue, Pittsburgh, PA 15213, USA}

\begin{abstract}
In this paper we consider strong Nash equilibria, in mixed strategies, for finite games. Any strong Nash equilibrium outcome is Pareto efficient for each coalition. First, we analyze the two--player setting. Our main result, in its simplest form, states that if a game has a strong Nash equilibrium with full support (that is, both players randomize among all pure strategies), then the game is strictly competitive. 
% NICOLA :: I changed the next sentence because, in the previous version of the paper, it was imprecise. 
This means that all the outcomes of the game are Pareto efficient and lie on a straight line with negative slope.
In order to get our result we use the indifference principle fulfilled by any Nash equilibrium,  and the classical KKT conditions (in the vector setting), that are necessary conditions for Pareto efficiency. 
%The flow doesn't work here because necessary conditions aren't suffient, and the conclusion here should be about being Pareto efficient??? 
%
%NICOLA :: it is not clear to me your comment. I added a sentence above to say that a strong Nash equilibrium is Pareto efficient (i.e., the second sentence of the abstract). This should provide the connection with "necessary condition for Pareto efficiency". Furthermore, another point could be that we use conditions that are only necessary and we provide a necessary--and--sufficient result. By the way, I don't see any contradiction. Indeed, indifference principle + KKT require the alignment of the outcomes (you can see this condition as a necessary condition), but in a strictly competitive game any NE is also an SNE (and therefore the alignment is also a sufficient condition).
Our characterization enables us to design a strong--Nash--equilibrium--finding algorithm with complexity in Smoothed--$\mathcal{P}$. So, this problem---that Conitzer and Sandholm [Conitzer, V., Sandholm, T., 2008. New complexity results about Nash equilibria. Games Econ. Behav. 63, 621--641] proved to be computationally hard in the worst case---is generically easy. Hence, although the worst case complexity of finding a strong Nash equilibrium is harder than that of finding a Nash equilibrium, once small perturbations are applied, finding a strong Nash is easier than finding a Nash equilibrium. Next we switch to the setting with more than two players. We demonstrate that a strong Nash equilibrium can exist in which an outcome that is strictly Pareto dominated by a Nash equilibrium occurs with positive probability.
Finally, 
we prove that games that have a strong Nash equilibrium where at least one player puts positive probability on at least two pure strategies are extremely rare: they are of zero measure.
\end{abstract}

%\begin{flushleft}
%{\sf 2010 Mathematics Subject Classification}: Primary 91A10, 91A05\\ Secondary 91A70.\\
%
%\vspace{.2cm} {\sf Keywords}:
%\end{flushleft}


\end{frontmatter}

\section{Introduction}
It is well known that in non--cooperative game theory~\cite{MSZ}, rationality implies that the players can be worse off  than they could be by collaborating. The most celebrated example is the prisoners' dilemma, where strictly dominant strategies for the players lead to a bad outcome for both. Aumann's \emph{strong Nash equilibrium}~\cite{A} gets around this paradox and also provides a solution concept that is robust against coalitional deviations.
A strategy profile is a Nash equilibrium~\cite{N} if no player has unilateral incentive to deviate, while a strategy profile is a strong Nash equilibrium if no coalition has incentive to deviate. It follows immediately that any strong Nash equilibrium outcome is weakly Pareto efficient for each coalition~\cite{P}. 

A simple refinement is  \emph{super strong Nash equilibrium}~\cite{R}, which requires strict Pareto efficiency for each coalition. There are classes of games that have a strong Nash equilibrium, but not a super strong Nash equilibrium~\cite{GM}, so the distinction between the two solution concepts is meaningful. 

The strong Nash equilibrium concept is commonly criticized as too demanding because it allows for unlimited private communication among the players, and in many games a strong Nash equilibrium does not exist. For these reasons among others, relaxations have been proposed.  A relaxation that we use in this paper is the concept of \emph{$k$--strong Nash equilibrium}. It is a Nash equilibrium where no coalition of $k$ or fewer agents has incentive to deviate~\cite{AFM}.  The rationale is that in many practical situations only small coalitions can form. Another relaxation is \emph{coalition--proof Nash equilibrium}, which is a Nash equilibrium that is resilient against those coalitional deviations that are self--enforcing~\cite{Bernheim19871}. Coalition--proof Nash equilibria are sometimes not Pareto efficient.

The property of having a strong Nash equilibrium is important for a game, but not ``usual''.  In this paper we consider the class of all games with a fixed number of players, each of them having a fixed, finite, number of strategies, and study  ``how many of the games'' have strong Nash equilibria and ``what is their geometry''. In the mathematical literature, ``how many of the games'' can be given various meanings~\cite{J}: there is the notion of ``zero measure''  set, or meager in the Baire sense, and also other, perhaps less known, ones: $\sigma$--porosity, sparseness and so on. Here we shall use the term ``negligible'', since, as we shall see, all these notions coincide in this setting. The reason for this is that the set of games that have a strong Nash equilibrium (in mixed strategies) can be characterized as a subset of a semialgebraic set, and for these sets all the above conditions to be ``small'' coincide. Dubey shows that strong Nash equilibria are generically in pure strategies in finite games~\cite{doi:10.1287/moor.11.1.1}. In the case of continuous convex games,  existence conditions have been provided~\cite{Nessah2014871}. In our paper, we extend the results provided by Dubey along several dimensions.\footnote{However, Dubey considers more general utility functions than we consider here.}

We take a different  approach than Dubey did. Ours is based on the application of the indifference principle and the Karush--Kuhn--Tucker conditions~\cite{M} that are necessary for weak Pareto efficiency. In terms of existence of strong Nash equilibrium, we provide an alternative proof to that of Dubey for the case with two players and we provide a stronger result than that provided by Dubey for the case with three or more players. Specifically, we show that even 2--strong Nash equilibria are generically in pure strategies and therefore $k$--strong Nash equilibria are generically pure for any $k\geq 2$. Our approach enables also the derivation of a number of new results.

Our main result is the precise description of the geometry of games that admit strong and super strong Nash equilibria, in the following sense. When  there are two players, the various outcomes of the games can be geometrically represented as points in the plane. We show that a strong Nash equilibrium in mixed strategies may exist only if all the outcomes of the game, restricted to the support of the equilibrium, lie on a straight line with non--strictly negative slope and therefore the game restricted to the support of the equilibrium must be either strictly competitive or have all the outcomes lying in a vertical or horizontal straight line.\footnote{Games in which all the outcomes lie in a vertical or horizontal straight line are degenerate games in which one player is indifferent over all the actions. That is, the player can be safely removed from the game.}  This also implies that in games admitting a strong Nash equilibrium all the outcomes of the game, restricted to the support of the equilibrium, are weakly Pareto efficient. Similar results hold for super strong Nash equilibrium. More precisely, a super strong Nash equilibrium in mixed strategies may exist only if all the outcomes of the game, restricted to the support of the equilibrium, lie on a straight line with \emph{strictly} negative slope and therefore the game restricted to the support of the equilibrium must be strictly competitive. Furthermore, in games admitting a super strong Nash equilibrium all the outcomes of the game, restricted to the support of the equilibrium, are \emph{strictly} Pareto efficient. We show instead that games with three or more players have different properties. Indeed, these games can have strong and super strong Nash equilibria in which the game restricted to the support of the equilibrium may contain outcomes that are strictly Pareto dominated.

We also provide results about the computational complexity of deciding whether a strong Nash equilibrium exists. Our geometric characterization of two--player games that admit strong and super strong Nash equilibria can be exploited to design an algorithm whose expected compute time, once a uniform perturbation $[-\sigma,\sigma]$ with $\sigma>0$ is applied independently to each entry of the bimatrix, is polynomial in the size of the game. Such an algorithm puts the problem of deciding whether a strong Nash equilibrium exists---and the problem of finding it---in Smoothed--$\mathcal{P}$ (the class of problems solvable in smoothed polynomial time), showing that these two problems are generically easy. 
A simple variation, omitted in this paper, applies to super strong Nash equilibrium showing that also the problems of deciding whether a super strong Nash equilibrium exists and of finding it are in Smoothed--$\mathcal{P}$. We recall that deciding whether there exists a strong Nash equilibrium is $\mathcal{NP}$--complete~\cite{sandholmComplexiy2008,aamasSNE2013} ($\mathcal{NP}$ is the class of non--deterministic polynomial time problems) and therefore no polynomial--time algorithm exists unless $\mathcal{NP}=\mathcal{P}$ ($\mathcal{P}$ is the class of polynomial time problems). Furthermore, finding a strong Nash equilibrium is harder than finding a Nash equilibrium. Indeed, finding a Nash equilibrium is $\mathcal{PPAD}$--complete~\cite{papastoc,Chen09:Settling} and $\mathcal{PPAD}\subset \mathcal{NP}$ unless $\mathcal{NP}=\text{co--}\mathcal{NP}$~\cite{Megiddo:1991:TFE:104800.104829} (co--$\mathcal{NP}$ is the class of complementary non--deterministic polynomial time), but it is commonly believed that $\mathcal{P}\subset \mathcal{PPAD}$  and therefore that no polynomial--time algorithm exists for finding a Nash equilibrium. Interestingly, finding a Nash equilibrium is not in Smoothed--$\mathcal{P}$ unless $\mathcal{PPAD}\subseteq \mathcal{RP}$~\cite{Chen09:Settling} ($\mathcal{RP}$ is the class of randomized polynomial time problems) and therefore it is commonly conjectured that finding a Nash equilibrium remains hard even as small perturbations are applied. Hence, although the worst case complexity of finding a strong Nash equilibrium is harder than that of finding a Nash equilibrium, once small perturbations are applied, finding a strong Nash is easier than finding a Nash equilibrium!

\section{Preliminaries}
First, we provide required notation on semialgebraic sets/functions that will be used in the rest of the paper.
\begin{defini}
A subset $A$ of an Euclidean space is called \rm algebraic \it if it can be described as a finite number of polynomial equalities. It is called \rm semialgebraic \it
if it can be described as a finite number of polynomial equalities and inequalities. A multivalued map between Euclidean spaces is called \rm algebraic (semialgebraic) \it if its graph is an algebraic (semialgebraic) set.
\end{defini}

For example, a circle in $\erre^2$ is an algebraic set,
%Do you mean just the perimeter of a circle??? ---> NICOLA :: yes, the perimeter, but "circle" should be correct. If you see "http://mathworld.wolfram.com/Circle.html", the circle is only the perimeter.
an interval in $\erre$ is a semialgebraic set. We shall use the idea of \it dimension \rm of a set, and this should be defined in general, but since here it is used for simplexes and/or affine spaces, it is enough to keep in mind that this coincides with the usual idea of dimension in linear analysis.

Here we need only two facts about semialgebraic maps:
\begin{itemize}
\item Given an algebraic set $A$ on $X\times Y$, its projection on each space $X$, $Y$ is semialgebraic (see \cite{C}, Section 1.3.1, the Tarski Seidenbreg theorem).
\item For any semialgebraic set--valued mapping $\Phi$ between two Euclidean spaces
$\Phi:E\rightrightarrows Y$, if  $\dim\Phi\left(x\right)\leq k$
for every $x\in E$, then $\dim\Phi\left(E\right)\leq\dim (E)+k$, where $\dim A$ denotes the dimension of a given set $A$ (see \cite{C}, Theorem 3.18).
\end{itemize}

We now introduce additional notation and definitions that we will use. 
For vectors $x,y$ in some Euclidean space, we use the notations $x\ge y$, $x>y$ and $x>>y$ to
say that $x_i\ge y_i$ for all $i$, $x_i\ge y_i$ for all $i$ and $x\ne y$, and $x_i> y_i$ for all $i$, respectively. 
We will use Pareto domination in the following setting.
We consider a function $F:\erre^k\to\erre^n$, and a set $U\subset \erre^k$.

\begin{defini}\label{(Strict-Pareto-dominance0)} A vector $\bar{x}=\left(\bar{x}_{1},\ldots,\bar{x}_{k}\right)\in U$
is {\rm weakly Pareto dominated} for the problem $(F,U)$ if there exists a vector $x\in U$ such that $$ F(x)> F(\bar x)$$
while $\bar{x}$
is said to be {\rm strictly Pareto dominated} for the problem $(F,U)$ if there exists $x\in U$ such that $$F(x)>>F(\bar x).$$
\end{defini}
On the basis of the concept of Pareto dominance,    Pareto efficiency can be defined.
\begin{defini}\label{(Strict-Pareto-Efficiency)} A vector $\bar{x}$
is {\rm strictly Pareto efficient} for the problem $(F,U)$  if there is no  $x$ that weakly Pareto dominates~$\bar{x}$, while it is
{\rm weakly Pareto efficient} for the problem $(F,U)$  if there is no  $x$  that strictly Pareto dominates~$\bar{x}$.
\end{defini}

We will be interested in the case where $U$ can be described in terms of affine inequalities. Thus, let us consider
$G:\erre^k\to\erre^l$, $H:\erre^k\to\erre^j$, $G,H$ linear, $b,c$ vectors of the right dimensions, and define $U$ to be
$$U=\{x\in\erre^k: G(x)\ge b,H(x)=c\}.$$
In such a case, if $\bar x\in U$ is weakly Pareto efficient, then the Karush--Kuhn--Tucker (KKT) conditions~\cite{M} state that
 there are vectors $\lambda,\mu,\nu$ that satisfy the following system:
\begin{subequations}\label{eq:CN_KKT}
\begin{eqnarray}
\textcolor{black}
\lambda^t\nabla F(x)+\mu^t\nabla G(x)+\nu^t\nabla H(x)&=&\un_n,\label{eq:KKT1} \\
\mu^{t}G(x) & = & 0,\label{eq:KKT2} \\
\mu & \geq & \un_{l},\\
\lambda & > & \un_j.
\end{eqnarray}
\end{subequations}

\medskip

Now we introduce the standard concepts from non-cooperative game theory that we will use in the paper.
\begin{defini}
(Strategic--form game) A \rm finite strategic--form game \it \cite{MSZ}
is a tuple $(N,A,U)$ where:
\begin{itemize}
\item $N=\left\{ 1,\ldots,n\right\} $ is the set of  players,
\item $A= A_{1}\times A_2\times \ldots \times A_{n}$ is the (finite) set of aggregate agents' actions:
$A_{i}$ is the set of actions available to agent $i$,
\item $U=\left\{ U_{1},\ldots,U_{n}\right\} $ is the set of agents' utility
tensors, where $U_{i}:A\to\erre$ is the utility function of agent $i$.
\end{itemize}
\end{defini}
We shall denote by $m_{i}\geq 1$
the number of actions in $A_{i}$, and by $a_{ij}$, $j\in\{1,...,m_{i}\}$,
a generic action;
$U$ is a $n\times m_{1}\times m_{2}\times...\times m_{n}$ tensor. A generic element of $U_i$ will be denoted by $U_i(i_1,\dots,i_n)$.

We denote  by
$\Delta_{i}$ the simplex of the mixed strategies over $A_{i}$, and by $x_{i}$ a generic mixed strategy of agent $i$: $x_i=\left(x_{i1},...,x_{im_{i}}\right)$.
 Given $x_{i}\in\Delta_{i}$ we denote
by $S_{i}(x_{i})$ its support, that is the set of actions
played with strictly positive probability, and, given a strategy profile $x$,  by $S(x)$ the support profile
$\left(S_{1}(x_1),\ldots,S_{n}(x_n)\right)$.

Given a strategy profile $x$, the utility of agent $i$ is
$$v_{i}(x)=\sum_{i_1,\dots,i_n}U_i(i_1,\dots,i_n)\cdot x_{i_1}\cdot \dots\cdot x_{i_n}:= x_i^tU_i\prod_{j\ne i}x_j.$$

Given a strategy profile $x$, we shall use $x_{-i}$ to denote the vector, with $n-1$ components,  $x_{-i}=(x_1,\dots,x_{i-1},x_{i+1},\dots,x_n)$ and we shall also write $x=(x_i,x_{-i})$.

\begin{defini}\label{(Nash-Equilibrium)} A strategy profile $\bar{x}=\left(\bar{x}_{1},\ldots,\bar{x}_{n}\right)$
is a \rm Nash equilibrium \it if, for each $i\in N$, $v_{i}\left(\bar{x}\right)\geq v_{i}\left(x_{i},\bar{x}_{-i}\right)$
for every $x_{i}\in$$\Delta_{i}$.
\end{defini}
More explicitly, $\bar{x}$ is a Nash equilibrium if $\bar{x}_{i}^{t}U_{i}\prod_{j\neq i}\bar{x}_{j}\geq x{}_{i}^{t}U_{i}\prod_{j\neq i}\bar{x}_{j}$
for every $i\in N$, for every $x_{i}\in\Delta_{i}$.

Given a mixed strategy profile $x$, and denoting by $U_i|_{S_i}$ the matrix containing only the rows in $S_i$, and by $U_i|_{S_i^c}$ its complement,
the problem of finding a Nash equilibrium can be expressed as the problem of finding
a profile strategy $x$ and, for all $i\in N$,  a real number $ v_{i}^{*}$ such that:
\begin{subequations}\label{eq:NE_cond}

\begin{eqnarray}
U_i|_{S_i}\prod_{j\neq i}x_{j}= & v_{i}^{*}\cdot\uun_{m_i} & \forall i\in N\label{eq:NE1}\\
U_i|_{S_i^c}\prod_{j\neq i}x_{j}\le & v_{i}^{*}\cdot\uun_{m_i} & \forall i\in N\label{eq:NE2}\\
x_{ij}\geq & 0 & \forall i\in N,\:\forall j\in\{1,..,m_{i}\}\label{eq:NE3}\\
x_{i}^{t}\cdot\uun_{m_i}= & 1 & \forall i\in N\label{eq:NE4}
\end{eqnarray}

\end{subequations}
\noindent
where $\uun_{m_i}$ is a column vector of $m_i$ positions with value 1 in every position. As is customary, we shall refer to the conditions in (\ref{eq:NE1})  as the \emph{indifference principle}.

\medskip

Finally, we are ready to define the key objects of study in this paper: Nash equilibria that are efficient for each coalition of players.
\begin{defini}\label{(Super-Strong-Nash)} A strategy profile $\bar{x}=\left(\bar{x}_{1},\ldots,\bar{x}_{n}\right)$ is a {\rm super strong Nash equilibrium} if it is a Nash equilibrium and it is strictly Pareto efficient for every coalition of players.
\end{defini}
\begin{defini}\label{(Strong-Nash)} A strategy profile $\bar{x}=\left(\bar{x}_{1},\ldots,\bar{x}_{n}\right)$ is a {\rm strong Nash equilibrium} if it is a Nash equilibrium and it is weakly Pareto efficient  for every coalition of players.
\end{defini}

We consider also relaxations of strong and super strong Nash equilibria in which we require the resilience against coalitions of size $k$ or less only, see \cite{AFM}.
\begin{defini}\label{(k-Super-Strong-Nash)} A strategy profile $\bar{x}=\left(\bar{x}_{1},\ldots,\bar{x}_{n}\right)$ is a {\rm k--super strong Nash equilibrium} if it is a Nash equilibrium and it is strictly Pareto efficient for every coalition of $k$ or fewer players, and it is 
 a {\rm k--strong Nash equilibrium} if it is a Nash equilibrium and it is weakly Pareto efficient  for every coalition of $k$ or fewer players.
\end{defini}

\section{Mixed strong Nash equilibria in two-player games}
In this section we consider the two--player setting. Later in the multi-player setting we will leverage some of these results. This is quite natural: since strong Nash equilibrium requires efficiency for all coalitions, we will use the results obtained in this section to the coalitions made by two players within the multi-player games.

We start the analysis focusing on a fully mixed strong Nash equilibrium. The game can be described by a bimatrix  \textcolor{black}{$(U_1,U_2^t)$ where $U_1=(u_1^{ij})$ and $U_2=(u_2^{ij})$}.

\subsection{Mixed strong Nash equilibria with full support}
Obviously, adding a constant to the payoffs of all players, and/or multiplying them by a positive constant, does not change the set of the strong Nash equilibria of the game. Thus, we shall assume without loss of generality that at a given equilibrium $x=(x_1,x_2)$ both players get zero: $x_1^tU_1x_2=0$ and $x_2^tU_2x_1=0$. \textcolor{black}{In other words, we are assuming that $v_i^* = 0$ in (\ref{eq:NE1}) and (\ref{eq:NE2}) for $i=1,2$.}

\medskip

\begin{prop}
Let $x$ be a fully mixed strong Nash equilibrium. Then it must fulfill the following system of linear equalities/inequalities, for some $\lambda=(\lambda_1,\lambda_2)$ and $\nu=(\nu_1,\nu_2)$:
\begin{subequations}
\begin{eqnarray}
\lambda_{2}\textcolor{black}{x_{2}^t}U_{2}+\nu_{1}\uun_{m_1} & = & \un_{\textcolor{black}{m_1}} \label{eq:1} \\
\lambda_{1}\textcolor{black}{x_{1}^t}U_{1} +\nu_{2}\uun_{m_2} & = & \un_{\textcolor{black}{m_2}} \label{eq:2}\\
\lambda >0 \label{eq:3}
\end{eqnarray}
\end{subequations}
\end{prop}
\pf \textcolor{black}{The proof follows from the application of the KKT conditions and the indifference principle. The elements used in the KKT conditions are
\begin{eqnarray*}
F(x) 	&	=	& 	\Big(	\begin{array}{cc}	f_1(x)	&	f_2(x)\end{array}	\Big) = \Big(	\begin{array}{cc}	x_1^tU_1x_2 	& 	x_2^tU_2x_1	\end{array}	 \Big)		\\
G(x) 	&	=	&	\Big(	\begin{array}{cc} x_{1}^t & x_{2}^t	\end{array}	\Big)																				 \\
H(x) 	&	= 	&	\Big(	\begin{array}{cc} x^t_{1}\uun_{m_1}-1 & x^t_{2}\uun_{m_2}-1 \end{array}	\Big)																 \\
\lambda&	=	&	\Big(	\begin{array}{cc} \lambda_1 & \lambda_2 \end{array}	\Big)																			 \\
\mu	&	=	&	\Big(	\begin{array}{cccccc} \mu_1 & 	\ldots	& \mu_{m_1}	&	\mu_{m_1+1}	&	\ldots	&	\mu_{m_1+m_2} \end{array}	\Big)					 \\
\nu	&	=	&	\Big(	\begin{array}{cc} \nu_1 & \nu_2 \end{array}	\Big)																			
\end{eqnarray*}
and therefore
\begin{eqnarray*}
\nabla F(x) 	&	=	& 	\left(	\begin{array}{cc}	U_1x_2 	& 	x_1^tU_1		\\ 	x_2^tU_2 	& 	U_2x_1		\end{array}	\right)	 \\
\nabla G(x) 	&	=	&	I_{m_1+m_2}																														 \\
\nabla H(x) 	&	=	& \left(	\begin{array}{cc}	\uun_{m_1}^t 	& 	\mathbf{0}_{m_2}^t		\\ 		\mathbf{0}_{m_1}^t 	& 	\mathbf{1}_{m_2}^t		 \end{array}	\right)\hspace{5.03cm}
\end{eqnarray*}
where $I_{m_1+m_2}$ is the identity matrix with $m_1+m_2$ rows and columns.}

From (\ref{eq:KKT1}) and from the fact that $x$ is fully mixed, we have $\mu = \un_{m_1+m_2}$. From the fact that $x$ is a Nash equilibrium and that, by assumption, $v_1^*=v_2^*=0$, we have $U_1x_2 = \un_{m_1}$ and $U_2x_1 = \un_{m_2}$. The claim follows straightforwardly. \newendproof

 \medskip

\begin{lemma}\label{villa}
\label{cap 5 lem:IP+KKT} Let $\Gamma:=(U_1,U_2^t)$ be a bimatrix game with $\max\{m_1,m_2\}\geq 2$.

Let  $x$ be a fully mixed super strong Nash equilibrium \textcolor{black}{of $\Gamma$} satisfying the system (\ref{eq:1}), (\ref{eq:2}), (\ref{eq:3}).  Then it satisfies the further conditions
$$U_{1}^{t}x_{1}=\un_{m_2},\qquad U_{2}^tx_2=\un_{m_1}$$
\textcolor{black}{also with strictly positive $\lambda_1,\lambda_2$.}

Let  $x$ be a fully mixed strong Nash equilibrium \textcolor{black}{of $\Gamma$} satisfying the system (\ref{eq:1}), (\ref{eq:2}), (\ref{eq:3}).  Then either it satisfies the further conditions
$$U_{1}^{t}x_{1}=\un_{m_2},\qquad U_{2}^tx_2=\un_{m_1}$$
or all entries of the bimatrix $(U_1,U_2^t)$ lie either on a vertical or on a horizontal  line through the origin.
%Somewhere in the paper it would be nice to have a figure (or otherwise explain) what "$(U_1,U_2^t)$ lie either on a vertical or on a horizontal  line through the origin" means???  ---> NICOLA :: I added a figure below
So, other than in the case where the entries of the bimatrix lie on a horizontal or vertical line, a fully mixed strong Nash  equilibrium satisfies the equations (\ref{eq:1}), (\ref{eq:2}),(\ref{eq:3}) also with strictly positive $\lambda_1,\lambda_2$.
\end{lemma}
\pf
We first prove the claim about super strong Nash equilibrium and then the claim about strong Nash equilibrium.

\emph{Super strong Nash equilibrium}. Suppose, without loss of generality, $\lambda_1>0$. We have the following steps.

Step 1. We show $x_1^tU_1=\mathbf{0}_{m_2}^t$. From (\ref{eq:2}) we know that $x_1^tU_1=\frac{\nu_2}{\lambda_1}\cdot\uun_{m_2}^t$. By assumption, we have $x_1^tU_1x_2 = 0$ and therefore $\frac{\nu_2}{\lambda_1} \uun^t_{m_2} x_2 =0$. Given that $\uun^t_{m_2} x_2= 1$, it follows that $\nu_2=0$ and thus $x_1^tU_1=\mathbf{0}_{m_2}^t$.

Step 2. We show $x_2^tU_2=\mathbf{0}_{m_1}^t$. Initially, we show $x_2^tU_2\leq\mathbf{0}_{m_1}^t$. Suppose for contradiction that there is a column~$j$ such that $x_2^tU_2^j>0$. Then, the strategy profile $y=(\bar{x}_1,x_2)$, where $\bar{x}_1$ is the pure $j$--th strategy of player~1, weakly Pareto dominates $x$. Indeed, $x_1^tU_1x_2 = x_2^tU_2x_1 = 0$ by assumption, whereas $\bar{x}_1^tU_1x_2 = 0$ and $x_2^tU_2\bar x_1>0$. As a result, we have a contradiction, $x$ being a super strong Nash equilibrium and therefore it cannot be weakly Pareto dominated. Thus, $x_2^tU_2\leq\mathbf{0}_{m_1}^t$. Finally, given that $x_2^tU_2x_1 = 0$ and $x_1\geq \mathbf{0}_{m_1}$, we have that  $x_2^tU_2\leq\mathbf{0}_{m_1}^t$ implies  $x_2^tU_2=\mathbf{0}_{m_1}^t$.

Step 3. Given that $x_2^tU_2=\mathbf{0}_{m_1}^t$, (\ref{eq:1}) is satisfied for $\nu_1 = 0$ and any $\lambda_2$, strictly positive values included. In other words, if $x$ a super strong Nash equilibrium, the system (\ref{eq:1}), (\ref{eq:2}),(\ref{eq:3}) admits a solution with strictly positive $\lambda_1,\lambda_2$.

\emph{Strong Nash equilibrium}. The above proof does not apply here. This is due to Step~2. It requires $x$ not to be weakly Pareto dominated, but a strong Nash equilibrium  may be weakly Pareto dominated. 
We observe in addition that, if the system (\ref{eq:1}), (\ref{eq:2}), (\ref{eq:3}) is satisfied at $x$ when $\lambda_1$ and $\lambda_2$ are  strictly positive, the proof easily follows from Steps~1--3. Indeed, Step~1 directly applies and we have $x_1^tU_1=\mathbf{0}_{m_2}^t$. In addition, by a simple variation of Step~1 when $\lambda_2 >0$ (instead of $\lambda_1 >0$), we have $x_2^tU_2=\mathbf{0}_{m_1}^t$. Therefore, we need to complete the proof for the case in which the system (\ref{eq:1}), (\ref{eq:2}), (\ref{eq:3}) is satisfied at $x$ only when one of the two components of $\lambda$, say $\lambda_2$, is vanishing, given that KKT conditions exclude the case in which both $\lambda_1$ and $\lambda_2$ are equal to zero. By applying Step~1, we already know that $x_1^tU_1=\mathbf{0}_{m_2}^t$, but, since $\lambda_2$ must be 0,  there must be $j$ such that $x_2^tU_2^j> 0$. We have the following steps.

Step 4. We show $u_1^{ji}=0$ for all $i$. Suppose for contradiction that there is $i$ such that $u_1^{ji}\neq 0$. Then, since $x_1^tU_1=\mathbf{0}_{m_2}^t$ and $x_1\geq 0$, there are $k,l$ such that $u_1^{jk}>0$ and $u_1^{jl}<0$. Since $x_2^tU_2^j> 0$, there is a neighborhood $B$ of $x_2$ (in the simplex) such that for all $y\in B$ it is $y^tU_2^j> 0$. Then it is possible to choose $\eps>0$ so small that $\bar x_2=(x_{21},\dots,x_{2k}+\eps,\dots,x_{2l}-\eps,\dots)$ is in $B$, but this is a contradiction, since $x$ is strictly Pareto dominated by $y=(\bar x_1,\bar x_2)$, where $\bar x_1$ is the pure $j$--th strategy of player~1. Indeed, $\bar x_2^t U_2 \bar x_1 > 0 = x_2^t U_2 x_1$ and $\bar x_1^t U_1 \bar x_2 > 0 = x_1^t U_1 x_2$. Therefore, it must be that $u_1^{ji}=0$ for all $i$, and this means that all entries of row $j$ of the bimatrix $(U_1,U_2^t)$ lie on a vertical line through the origin.

Step 5. We show $u_1^{il}=0$ for all $i$ and $l$. Suppose for contradiction that $u_1^{il}>0$ (observe, $i\ne j$). Then consider the following strategy $\bar x_1$ for player~1:
$$
\begin{cases}
\bar x_{1k}&=0 \qquad\qquad \text{if } k\ne i,j\\
\bar x_{1i}&=\sigma>0 \\
\bar x_{1j}&=1-\sigma
\end{cases}
$$
with $\sigma$ so small that $(1-\sigma)\sum_k x_{2k}u_2^{kj}+\sigma \sum_k x_{2k}u_2^{ki}>0$. Consider also the strategy $\bar x_2$ for player~2:
$$
\begin{cases}
\bar x_{2k}&=x_{2k} \qquad\qquad \text{if } k\ne l,m\\
\bar x_{2l}&=x_{2l}+\eps \\
\bar x_{2m}&=x_{2m}-\eps
\end{cases}
$$
for $m$ such that $u_1^{im}<0$ \footnote{Remember: $x_2^tU_1^i=0$ and thus, since $u_1^{il}>0$ there must be also a negative entry in the row $i$.} and $\eps$ so small that $\bar x_2\in B$ and $(1-\sigma)\sum_k \bar x_{2k} u_2^{kj}+\sigma \sum_k \bar x_{2k} u_2^{ki}>0$. By construction, we have $\bar x_2^t U_2 \bar x_1 > 0$. In addition, we have $\bar x_1^t U_1 \bar x_2 = \bar x_1^t U_1 x_2 + (1-\sigma)(\epsilon u_1^{jl} - \epsilon u_1^{jm}) + \sigma (\epsilon u_1^{il} - \epsilon u_1^{im})=(1-\sigma)(\epsilon u_1^{jl} - \epsilon u_1^{jm}) + \sigma (\epsilon u_1^{il} - \epsilon u_1^{im})$, but,  since $u_1^{jl}=u_1^{jm}=0$ as showed in Step~4 and $u_1^{il}>0$ and $u_1^{im}<0$ by assumption, $\bar x_1^t U_1 \bar x_2>0$ and therefore we have a contradiction, given that strategy profile $\bar x$ strictly Pareto dominates~$x$. Thus player one gets zero at every outcome. This means that all the outcomes  of the bimatrix lie on the vertical axis.\newendproof

 \medskip

We can now prove the main theorem of this section.

\begin{thm} \label{oggi} Let $\Gamma:=(U_1,U_2^t)$ be a bimatrix game  with $\max\{m_1,m_2\}\geq 2$.

Let $x$ be a fully mixed super strong Nash equilibrium of $\Gamma$, providing zero utility to both players. Then all the entries $u^{ij}:=(u_1^{ij}, u_2^{ji})$ of the bimatrix lie on the same straight line with strictly negative slope, passing through the origin.

Let $x$ be a fully mixed strong Nash equilibrium of $\Gamma$,  providing zero utility to both players. Then all the entries $u^{ij}:=(u_1^{ij}, u_2^{ji})$ of the bimatrix lie on the same straight line with non--strictly positive slope, passing through the origin.
\end{thm}
\pf
We first prove the claim about super strong Nash equilibrium and then the claim about strong Nash equilibrium.

\emph{Super strong Nash equilibrium}. We have the following steps.

Step 1. We show that all the entries $u^{ij}$ on a fixed row~$i$ lie on the same straight line through the origin and the same holds for a fixed column~$j$ and these straight lines cannot be vertical or horizontal. Let us focus on the rows, the same reasoning can be applied to the columns. Suppose for contradiction that for a given row~$i$ the entries $u^{ij}$ do not lie on the same straight line. Consider the convex hull of the entries $u^{ij}$ in the space of the players' expected utilities. This is a polygon whose vertices are a subset of the entries $u^{ij}$. For each point $(v_1,v_2)$ of the polygon, there is a strategy $\bar{x}_2$ such that $\bar{x}_1^tU_1\bar{x}_2=v_1$ and $\bar{x}_2^tU_2\bar{x}_1=v_2$, where $\bar{x}_1$ is the $i$--th pure strategy of player~$1$. In addition, by Lemma~\ref{villa}, we know $U_1x_2 = U_2^t x_2 =\un_{m_1}$ and therefore the point $(0,0)$, corresponding to the players' expected utilities given by strategy profile~$x=(x_1,x_2)$, is in the interior of the polygon given that $x_2$ is fully mixed. Thus, we have a contradiction since there must be a strategy profile $\bar{x}=(\bar{x}_1,\bar{x}_2)$ such that $\bar{x}_1^tU_1\bar{x}_2 = \epsilon_1>0$ and $\bar{x}_2^tU_2\bar{x}_1 = \epsilon_2> 0$ for some small $\epsilon_1,\epsilon_2$, where $(\epsilon_1,\epsilon_2)$ is on the boundary of the polygon, and therefore $x$ is strictly Pareto dominated by $\bar{x}$. This shows that all the entries $u^{ij}$ for a given row~$i$ lie on the same straight line passing through the origin. In a similar way, we obtain that all the entries $u^{ij}$ for a given column~$j$ lie on the same straight line passing through the origin

Finally, each straight line has strictly negative slopes or collapses to a single point. Otherwise, if there is at least a straight line with strictly positive slope, there must be a strategy profile $\bar{x}=(\bar{x}_1,\bar{x}_2)$, as defined above, that weakly Pareto dominates~$x$.

Step 2. We show that all the entries $u^{ij}$ lie on the same straight line. Suppose for contradiction, without loss of generality, that the entries of two rows $i,j$ belong to different lines passing through the origin. Then, there must be two columns $h,k$ with $h\neq k$ such that the points $u^{ih}$ and $u^{jk}$ do not belong to the same straight line passing trough the origin. Note that $h$ must be different than $k$ due to Step~1 that forces all the entries of the same column to be on the same straight line passing through the origin. Let us now consider the sub--bimatrix
\[
\begin{pmatrix}
u^{ih} 	& 	u^{ik} 	\\
u^{jh} 	& 	u^{jk}
\end{pmatrix}.
\]
We assume that the segment joining $u^{ih}$ and $u^{jk}$ intersects the first orthant. 
%What does this mean?  What are the orthants and how do you number them? Would be nice to have a figure??? ---> NICOLA :: orthant means, in the case of 2 dimensions, quadrant. I added a figure below that shows the numbering of the orthant
This assumption is without loss of generality, because, for each pair of $u^{ih}$ and $u^{jk}$ such that the segment joining them does intersect the third orthant there is a pair $u^{ih'}$ and $u^{jk'}$ such that the segment joining them does intersect the first orthant.

We must have $u^{ik}= (0,0)$. If $u^{ik}\neq (0,0)$, then $u^{ih}$ and $u^{jk}$ would be on the same straight line passing trough $u^{ih}$ and $(0,0)$; this contradicts the assumption that $u^{ih}$ and $u^{jk}$ do not belong to the same straight line passing through the origin. The same reasoning applies to $u^{jh}$, and therefore $u^{jh}= (0,0)$. Thus, what remains to be done is to consider the case
\[
\begin{pmatrix}
u^{ih} 	& 	(0,0)	 	\\
(0,0	) 	& 	u^{jk}
\end{pmatrix}.
\]
What we need to prove is that this configuration leads to a contradiction. Since we assumed that the segment joining $u^{ih}$ and $u^{jk}$ intersects the first orthant,
there is $t>0$  such that
$$t u^{ih}+ (1-t)u^{jk}=(2\epsilon,2\epsilon)$$ for some $\epsilon>0$. Consider the strategy profile $\bar x=(\bar{x}_1,\bar{x}_2)=[(t,1-t),(\frac{1}{2}, \frac{1}{2})]$. Then $\bar{x}_1^t\bar U_i \bar{x}_1=\epsilon>0$ for $i=1,2$, contradicting the fact that $x$ is a super strong Nash equilibrium.

\emph{Strong Nash equilibrium}. In this case, the above Step~1 applies with the exception that the straight lines can be vertical or horizontal as stated by Lemma~\ref{villa}, while Step~2 applies here without exception. This concludes the proof.\newendproof

\medskip

\begin{cor}\label{oggi1} Let $\Gamma:=(U_1,U_2^t)$ be a bimatrix game  with $\max\{m_1,m_2\}\geq 2$.

Suppose $\Gamma$ has a super strong Nash equilibrium with fully mixed strategies, then it is a strictly competitive game.

Suppose $\Gamma$ has a strong Nash equilibrium with fully mixed strategies, then either it is a strictly competitive game or all entries of the bimatrix $(U_1,U_2^t)$ lie either on a vertical or on a horizontal  line through the origin.
\end{cor}

\medskip

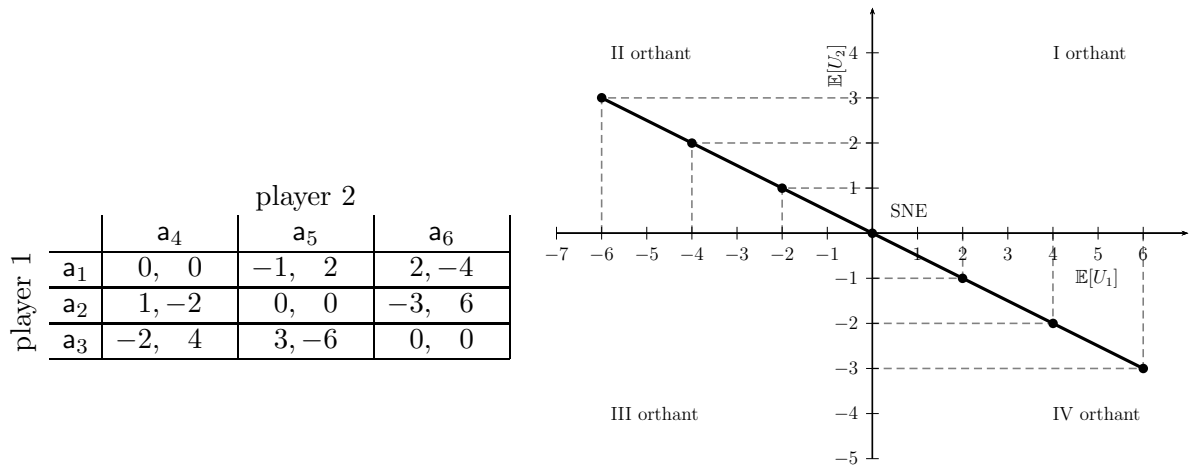
\begin{figure}[h]
\begin{minipage}{4.5cm}
\[
\begin{small}
\begin{array}{rr|c|c|c|}
\multicolumn{2}{c}{}	&	\multicolumn{3}{c}{\textnormal{player 2}} \\
	&		&	\mathsf{a}_4	&	\mathsf{a}_5	&	\mathsf{a}_6	\\ \cline{2-5}
\multirow{3}{*}{\begin{sideways}player 1\end{sideways}}	&	\mathsf{a}_1	&	\textcolor{white}{-}0,\textcolor{white}{-}0\textcolor{white}{-} 	& 	-1,\textcolor{white}{-}2\textcolor{white}{-}	&	\textcolor{white}{-}2,-4\textcolor{white}{-}	\\ \cline{2-5}
	&	\mathsf{a}_2	&	\textcolor{white}{-}1,-2\textcolor{white}{-}	&	\textcolor{white}{-}0,\textcolor{white}{-}0\textcolor{white}{-}	&	-3,\textcolor{white}{-}6\textcolor{white}{-}	\\ \cline{2-5}
	&	\mathsf{a}_3	&	-2,\textcolor{white}{-}4\textcolor{white}{-}	&	\textcolor{white}{-}3,-6\textcolor{white}{-}	&	\textcolor{white}{-}0,\textcolor{white}{-}0\textcolor{white}{-}	\\ \cline{2-5}
\end{array}
\end{small}
\]
\end{minipage}
\begin{minipage}{5cm}
\begin{pspicture}*(-7,-3.5)(7,3)
\scalebox{0.6}{

\psline[linecolor=gray,linewidth=1pt,linestyle=dashed]{-}(-6,0)(-6,3)(0,3)
\psline[linecolor=gray,linewidth=1pt,linestyle=dashed]{-}(6,0)(6,-3)(0,-3)
\psline[linecolor=gray,linewidth=1pt,linestyle=dashed]{-}(-4,0)(-4,2)(0,2)
\psline[linecolor=gray,linewidth=1pt,linestyle=dashed]{-}(4,0)(4,-2)(0,-2)
\psline[linecolor=gray,linewidth=1pt,linestyle=dashed]{-}(-2,0)(-2,1)(0,1)
\psline[linecolor=gray,linewidth=1pt,linestyle=dashed]{-}(2,0)(2,-1)(0,-1)

\psaxes[]{->}(0,0)(-7,-5)(7,5)

\psline[linecolor=black,linewidth=2pt]{-}(-6,3)(6,-3)

\psdots[linecolor=black,dotsize=6pt](-6,3)(6,-3)(-4,2)(4,-2)(-2,1)(2,-1)(0,0)

\uput{0}[0](0.4,0.5){SNE}
\uput{0}[0](4.5,-1.0){$\mathbb{E}[U_1]$}
\rput[tr]{90}(-1.0,4.2){$\mathbb{E}[U_2]$}

\uput{0}[0](4,4){I orthant}
\uput{0}[0](-5.8,4){II orthant}
\uput{0}[0](-5.8,-4){III orthant}
\uput{0}[0](4,-4){IV orthant}
}
\end{pspicture}
\end{minipage}
\caption{Example of a 2--player strictly competitive game (left) and its Pareto curve (right) in which the expected utilities of the players at the SNE are $(0,0)$.}
\label{exam1}
\end{figure}

A strictly competitive game is a game such that, up to a suitable rescaling of the utility of one player, the game is zero sum~\cite{ADP}. An example of strictly competitive game is depicted in Fig.~\ref{exam1}. \textcolor{black}{Corollary~\ref{oggi1} implies that a necessary condition for the existence of a super strong Nash equilibrium in fully mixed strategies is that all the outcomes are \emph{strictly} Pareto efficient and that a necessary condition for the existence of a strong Nash equilibrium in fully mixed strategies is that all the outcomes are \emph{weakly} Pareto efficient. Furthermore, Corollary~\ref{oggi1} requires also that all the outcomes lie on the same straight line. This is  the simplest case in which all the outcomes are Pareto efficient and, in this case, any Nash equilibrium is also strong.}

\subsection{Mixed strong Nash equilibria without full support}
The results obtained above can be rephrased also in the case of Nash equilibria with no full support. This is quite simple to do, for the following reason. Given the game $\Gamma:=(U_1,U_2^t)$, if we consider a mixed strategy (super) strong Nash equilibrium profile $x$, then it is clear that its restriction on the support $S(x)$ must necessarily be a (super) strong Nash equilibrium for the game $\Gamma_r$ whose bimatrix is the restriction of $(U_1,U_2^t)$ to $S(x)$. Thus the above result applies, and \textcolor{black}{all the outcomes of the restricted game must lie on straight line with strictly negative slope, in the case of super strong Nash, or on a straight line with non--strictly positive slope, in the case of strong Nash}. Observe also that the games having a Nash equilibrium with one player using a pure strategy and the other one a mixed, non pure, strategy, is itself zero measure, without requiring any form of efficiency of the Nash equilibrium. In particular, this easily implies the following existence result.
\begin{thm}
In the space of all the bimatrix games of fixed dimension, the set of games having a strong Nash equilibrium in which at least one player plays a mixed strategy has zero measure.
\end{thm}

More precisely, we have shown that this set of games with strong Nash equilibria where at least one player uses a mixed strategy is contained in a subspace of dimension strictly less than the whole space. The above theorem obviously implies the same result also for super strong Nash equilibrium, given that super strong Nash equilibrium is a refinement of strong Nash equilibrium. We recall that the set of games with at least one strong Nash equilibrium contains an open set and therefore games lying in the interior of such set do necessarily possess a strong Nash equilibrium in pure strategies.

\subsection{Smoothed--$\mathcal{P}$ complexity and strong Nash equilibrium}

The characterization provided by Theorem~\ref{oggi} suggests a simple algorithm to find an SNE in bimatrix games. Initially, we introduce two conditions the algorithm exploits:
\begin{description}
\item[$Condition~1:=$] ``there is 2x2 sub bimatrix of $(U_1,U_2^t)$ in which all the entries lie on a line'';
\item[$Condition~2:=$] ``there is 2x1 sub bimatrix of $(U_1,U_2^t)$ in which all the entries lie on a vertical line or there is 1x2 sub bimatrix of $(U_1,U_2^t)$ in which all the entries lie on a horizontal line''.
\end{description}

\begin{algorithm}[h]\caption{SNE--finding($U_1,U_2^t$)}\label{alg:SNE}
\begin{algorithmic}[1]

\ForAll{pure--strategy profiles $x$}
\If{$x$ is a Nash equilibrium}
\If{$x$ is Pareto efficient}
\State \textbf{return} $x$
\EndIf
\EndIf
\EndFor
\If{$Condition~1$ holds or $Condition~2$ holds}
\ForAll{support profiles $\bar{S}$}
\If{there is a Nash equilibrium $x^*$ with $S(x)=\bar{S}$ (in case of multiple equilibria take $x^*$ as the equilibrium maximizing the social welfare)}
\If{$x^*$ is Pareto efficient}
\State \textbf{return} $x^*$
\EndIf
\EndIf
\EndFor
\EndIf
\State \textbf{return} NonExistence
\end{algorithmic}
\end{algorithm}

We recall that, given a strategy profile $x$, it can be verified whether $x$ is weak Pareto efficient (Steps~3 and~11) by means of the algorithm described in~\cite{aamasSNE2013}, while, given a support profile $\bar{S}$, it can be verified whether there is a Nash equilibrium $x^*$ with $S(x)=\bar{S}$ (Step~10) by means of linear programming as shown in~\cite{DBLP:journals/geb/PorterNS08}.

\begin{thm} \label{smoothed} Let $\Gamma:=(U_1,U_2^t)$ be a bimatrix game. The problem of finding a strong Nash equilibrium of $\Gamma$ is in Smoothed--$\mathcal{P}$.
\end{thm}
\pf
In order to show that the problem of finding a strong Nash equilibrium of $\Gamma$ is in Smoothed--$\mathcal{P}$, we need to show that the expected time of Algorithm~\ref{alg:SNE} once a perturbation $[-\sigma,+\sigma]$, where $\sigma>0$, with uniform probability is applied  to each entry of the bimatrix independently is polynomial in the size of the game (i.e., $m_1$ and $m_2$)~\cite{DBLP:journals/corr/abs-1202-1936}.

We initially observe that Steps~1--7 of Algorithm~\ref{alg:SNE} have complexity polynomial in the size of the game, the number of pure--strategy profiles being $m_1\cdot m_2$ and Steps~2--3 requiring polynomial time in $m_1$ and $m_2$. Then, we observe that, once a perturbation $[-\sigma,+\sigma]$ with uniform probability is applied  to each entry of the bimatrix independently, $Condition~1$ and $Condition~2$ are verified with zero probability. Then, Steps~9--15, that require exponential time in $m_1$ and $m_2$ in the worst case, are executed with zero probability. This proves that the expected time of Algorithm~\ref{alg:SNE} is polynomial in $m_1$ and $m_2$.
\newendproof

This shows that the problems of deciding whether a strong Nash equilibrium exists and of finding it are generically easy. A simple variation of Algorithm~\ref{alg:SNE} can be designed to show that also the problems of deciding whether a super strong Nash equilibrium exists and of finding it are generically easy.

\section{Setting with more than two players}

In this section, we extend our analysis to a generic $n$--player game, in order to see what conditions concerning the existence of a mixed--strategy strong Nash  we can obtain. The application of the indifference principle and of the KKT conditions leads to an equation system that, differently from the two--player case, is not linear. For this reason, we resort to semi--algebraic set--valued mappings. Here the result.

\begin{thm}\label{threeormoreplayers}
In the space of all $n$--player games of fixed dimension, the set of games having a 2--strong Nash equilibrium in which at least one player plays a mixed strategy has zero measure.
\end{thm}
\pf
We initially prove the theorem for super strong Nash equilibrium and subsequently for strong Nash equilibrium.

\emph{Super strong Nash equilibrium}. Consider a generic support profile in which at least one player randomizes over at least two actions. Let $s_i$ the number of actions in the support of player~$i$ and, without loss of generality, let $s_n=\max\limits_{i\in N}\{s_i\}\geq 2$.

By applying the IP, we derive the following equations:
\begin{equation}\label{IPmorethreeplayers}
\begin{array}{lllll}
U_i \prod\limits_{j\neq i}x_j 	&	= 	&	v_i^* \mathbf{1}_{s_i}	&	\forall i \in N
\end{array}
\end{equation}

By applying the KKT conditions to all the coalitions of two players of the form $(i,n)$ for all $i\leq n-1$, we derive the following equations:
\begin{equation}\label{Eccoci}
\begin{array}{llllll}
\lambda_{i,(i,n)} v_i^* \mathbf{1}_{s_i} 				&	+	&	\lambda_{n,(i,n)} U_n \prod\limits_{j\neq i}x_j 			& 		= 	&	\nu_{i,(i,n)} \mathbf{1}_{s_i}		 &	\forall i\in N\setminus\{n\}	
\end{array}
\end{equation}

\begin{equation}\label{eq:Eccocibis}
\begin{array}{llllll}
	\lambda_{i,(i,n)} U_i \prod \limits_{j\neq n}x_j 			&	+  	&	\lambda_{n,(i,n)} v_n^* \mathbf{1}_{s_n}				& 		= 	&	\nu_{n,(i,n)} \mathbf{1}_{s_n}		 &	\forall i\in N\setminus\{n\}		
\end{array}
\end{equation}

From equations (\ref{eq:Eccocibis}), we derive the following equations, for some $v_i^*$:
\begin{equation}\label{PEmorethreeplayers}
\begin{array}{lllll}
U_i \prod\limits_{j\neq n}x_j 	&	= 	&	v_i^* \mathbf{1}_{s_n}	&	\forall i \in N\setminus\{n\}		
\end{array}
\end{equation}

From the sets of Equations ~(\ref{IPmorethreeplayers}) (used for $i=n$) and~(\ref{PEmorethreeplayers}) (used for all $i\neq n$), we extract the following subset of equations:
\begin{equation}\label{equationsystem8}
\begin{array}{cccc}
U_1 \prod\limits_{j\neq n}x_j 	&	= 	&	v_1^* \mathbf{1}_{s_n}	\\
U_2 \prod\limits_{j\neq n}x_j 	&	= 	&	v_2^* \mathbf{1}_{s_n}	\\
\vdots					&	\vdots	&	\vdots				\\
U_n \prod\limits_{j\neq n}x_j 	&	= 	&	v_n^* \mathbf{1}_{s_n}	\\
\end{array}
\end{equation}

The above set is composed of $ns_n$ equations. Call $k$ the row describing the $k$--th linear equation ($k=1,\ldots,ns_n$), where $U_1, U_2,\ldots,U_n$ are unknown, while $x_i$ and $v_i^*$ are fixed for every $i \in N$.

We denote the space of the $s_{1}\times s_{2}\times \ldots s_{n-1}$ tensors by $\mathbf{M}^{s_1\times s_2\times \ldots \times s_{n-1}}$. Consider the algebraic set--valued mapping
\[
\Phi:\Delta^{s_1}\times \Delta^{s_2}\times \ldots \times \Delta^{s_{n-1}}\times \mathbb{R}^n \rightrightarrows (\mathbf{M}^{s_1\times s_2\times \ldots \times s_{n-1}})^{n s_n}
\]
defined as
\[
\Phi(x_1,x_2,\ldots,x_n,v_1^*,v_2^*,\ldots,v_n^*)= \left\{(B_1,B_2,\ldots,B_{ns_n}): B_k \prod\limits_{j\neq n}x_j - v^*_{\left\lceil	\frac{k}{s_n}	 \right\rceil}=0, \forall k \leq ns_n \right\}.
\]
We want to show that $\Phi(\Delta^{s_1}\times \Delta^{s_2}\times \ldots \times \Delta^{s-1}\times \mathbb{R}^n)$ is negligible.

Observe that, for non null $\left(x_{1},x_{2},\ldots,x_{n-1}\right)$ to ensure
\[
\left(B_{1},B_{2},\ldots,B_{ns_{n}}\right)\in\Phi\left(x_1,x_2,\ldots,x_n,v_1^*,v_2^*,\ldots,v_n^*\right)
\]
simply requires each tensor $B_{i}$, independently, to lie in some hyperplane in a space of $ns_1s_2\ldots s_n$ dimensions, and hence
\[
\dim\Phi(x_1,x_2,\ldots,x_n,v_1^*,v_2^*,\ldots,v_n^*)=n s_1s_2\ldots s_n - ns_n.
\]

From this we get that
\begin{align*}
\dim\Phi(\Delta^{s_1}\times  \ldots \times \Delta^{s_{n-1}}\times \mathbb{R}^n)		& 	=            \\
                                                                            & = n s_1\ldots s_n - ns_n + (s_1-1) + \ldots + (s_{n-1}-1)+n 		\\
																			&	\leq n s_1\ldots s_n -ns_n +(s_1+\ldots s_{n-1})+1
                                                                            \\ & \le 	n s_1\ldots s_n -ns_n+(n-1)s_n+1						\\
																			&	< n s_1\ldots s_n
\end{align*}
since $s_n\geq 2$. Therefore the set of tensors satisfying the conditions is negligible with respect to the space of all the utility tensors. Given that the above property holds for every (non--pure) support profile and given that the support profiles are finite, the thesis of the theorem follows straightforwardly.

\emph{Strong Nash equilibrium}. The above proof for super strong Nash equilibrium applies directly also for strong Nash equilibrium whenever multipliers $\lambda_{i,(i,n)}$ for every $i$ are strictly positive. Indeed, in this case, for each coalition $(i,n)$ we have $U_i \prod\limits_{j\neq n} x_j= v^*_i \mathbf{1}_{s_n}$ and therefore we can derive entirely the equation system~(\ref{equationsystem8}). However, we know that a strong Nash equilibrium may satisfy KKT conditions even when some multiplier is zero and such a null multiplier may be, in principle, $\lambda_{i,(i,n)}$ for some coalition $(i,n)$. In this case, the above proof for super strong Nash equilibrium does not apply directly for strong Nash equilibrium. We show below how the proof can be modified to capture this case.

Initially, consider a generic support profile in which all the  players randomize over at least two actions. Let $s_i$ the number of actions in the support of player~$i$ and, without loss of generality, let $s_n=\max\limits_{i\in N}\{s_i\}\geq 2$. We consider all the coalitions of two players in which one player is~$n$. Suppose that, for at least one coalition $(i,n)$, KKT conditions are satisfied only for $\lambda_{i,(i,n)}=0$. The proof of Lemma~\ref{cap 5 lem:IP+KKT} shows that in this case $U_n \prod\limits_{j\neq i,n}x_j= v^*_n \mathbf{M}_1^{s_i\times s_n}$ where $\mathbf{M}_1^{s_i\times s_n}$ is a $s_i\times s_n$ matrix of ones. Indeed, fixed $x_j$ for every $j\neq i,n$, the game reduces to a bimatrix game $(U_i', U'_n )$ between players $i$ and $n$ where $U_i' =U_i  \prod\limits_{j\neq i,n}x_j$ and $U_n' =U_n \prod\limits_{j\neq i,n}x_j$ and Lemma~\ref{cap 5 lem:IP+KKT} shows that all the entries of $U'_n$ are equal to $v^*_n$ when $\lambda_i=0$. We use these conditions together with the equations due to the indifference principle for all the players $j\neq n$, obtaining:
\begin{align}
U_j \prod\limits_{k\neq j}x_k 	&	= 	&&	v_j^* \mathbf{1}_{s_j}					&	\forall j \in N\setminus\{n\}		\label{tempstrongproof1}	\\
U_n \prod\limits_{k\neq i,n}x_k 	&	= 	&&	v_n^* \mathbf{M}_{1}^{s_i\times s_n}		&							\label{tempstrongproof2}
\end{align}

As in the proof for the super strong Nash equilibrium, call $\Phi$ the algebraic set--valued mapping that, given $(x_1,x_2,\ldots,x_n,v_1^*,v_2^*,\ldots,v_n^*)$ returns the set of tensors $U_1,\ldots,U_n$ (each of size $\prod_j s_j$) satisfying the above set of equations (\ref{tempstrongproof1})--(\ref{tempstrongproof2}). It can be seen that $\dim \Phi(\Delta^{s_1}\times \Delta^{s_2}\times \ldots \times \Delta^{s_{n}}\times \mathbb{R}^n) = n\prod_{j}s_j + s_n(1-s_i)<n\prod_{j}s_j$ and therefore $\Phi$ is negligible.

Finally, consider a generic support profile in which at least one player randomizes over at least two actions. Extract the reduced game containing only the actions belonging to the supports of the players and in which all the non--randomizing players have been discarded. Then, apply the above arguments to the reduced game in which all the players randomize. We obtain that the reduced game cannot be generic and therefore the original game cannot be generic.

Given that the above property holds for every (non--pure) support profile and given that the support profiles are finite, the thesis of the theorem follows straightforwardly.
\newendproof

\medskip
Let us note that Dubey shows in~\cite{doi:10.1287/moor.11.1.1} that, generically, in a Pareto efficient Nash equilibrium at least one player plays a pure strategy and that in a strong Nash equilibrium all the players play pure strategies. The above result  shows that generically pure--strategy equilibria are the only possible ones even when the resilience to multilateral deviations is required only for coalitions of two or less players.

To conclude, let us make the following observation. In the two--player case, our result shows that existence of a mixed strong Nash equilibrium requires  that all the outcomes restricted to the support of the equilibrium are weakly Pareto efficient. Interestingly, this is no longer true with three players.
\begin{prop}\label{prop:dominatedoutcomes3players}
Mixed strong Nash equilibria of three--player games may have outcomes that are strictly Pareto dominated. The same holds for super strong Nash equilibrium.
\end{prop}
\noindent \textbf{Proof sketch}
We consider the following game with three players, in which every player has two available
actions:
\[
M_1=\begin{pmatrix}
(2,0,0) & (0,2,0) \\ (0,0,2) & (0,0,0)\\
\end{pmatrix},
M_2=\begin{pmatrix}
(0,0,0) & (0,0,2) \\ (0,2,0) & (2,0,0)\\
\end{pmatrix}.
\]

The profile strategy $\left(\left(\frac{1}{2},\frac{1}{2}\right),\left(\frac{1}{2},\frac{1}{2}\right),\left(\frac{1}{2},\frac{1}{2}\right)\right)$
is a super strong Nash equilibrium. For, it is easy to see that
$\left(\left(\frac{1}{2},\frac{1}{2}\right),\left(\frac{1}{2},\frac{1}{2}\right),\left(\frac{1}{2},\frac{1}{2}\right)\right)$
is a Nash equilibrium, with value $v^{*}=\left(\frac{1}{2},\frac{1}{2},\frac{1}{2}\right)$. \textcolor{black}{Then, it is possible to prove that neither the  coalitions made by two players nor the grand coalition have incentive to deviate from the Nash equilibrium. Calculations are straightforward, but long, and thus we report them in~\ref{appendix}. Therefore, $\left(\left(\frac{1}{2},\frac{1}{2}\right),\left(\frac{1}{2},\frac{1}{2}\right),\left(\frac{1}{2},\frac{1}{2}\right)\right)$ is a super strong Nash equilibrium.}

\textcolor{black}{Finally, it is easy to see that the strong Nash equilibrium strictly Pareto dominates outcome $(0,0,0)$. This concludes the proof.}
\newendproof

\medskip
 Observe that in the example above there are repeated entries in the matrices, but small perturbations of the two matrices
provide the same results, having all triples of the two matrices different each other. Thus the above example is generic.
Proposition \ref{prop:dominatedoutcomes3players} shows that the geometric characterization in terms of alignment of the bimatrix entries of two--player games admitting (super) strong Nash equilibria cannot be extended to the case of three or more players. Furthermore, it is not clear whether there is a simple geometric characterization in the case of three or more players. This leaves the question of setting the problem of deciding whether there is a (super) strong Nash equilibrium in Smoothed--$\mathcal{P}$ open.

\section{Concluding remarks}
In this paper we analyzed the problem of characterizing the set of finite games having strong Nash equilibria, extending the results provided by Dubey in~\cite{doi:10.1287/moor.11.1.1}. Our main result concerns the characterization of two--player games admitting strong and super strong Nash equilibria. We showed that the game restricted to the support of the equilibrium must be strictly competitive in the case of super strong Nash equilibrium and must be either strictly competitive or with all the outcomes that lie on a horizontal or vertical line in the case of strong Nash equilibria. This implies that all the outcomes of the game restricted to the support of the equilibrium must be Pareto efficient. This is no longer true in the case of three or more players, where instead the support of even super strong Nash equilibria may contain strictly Pareto dominated outcomes. For these games, we show that even 2--strong Nash equilibria are generically in pure strategies. Our geometric characterization of two--player games admitting a strong Nash equilibrium leads to the design of a simple algorithm that puts the problems of deciding whether there is a strong Nash equilibrium and of finding it in Smoothed--$\mathcal{P}$. This shows that such problems are generically easy.

The main question we leave open concerns the characterization of the geometry of games with three or more players admitting strong Nash equilibrium. It is not clear whether or not their geometry can lead to the design of  a Smoothed--$\mathcal{P}$ algorithm.

\bigskip

\noindent {\bf Acknowledgements.} The authors gratefully acknowledge Prof. A. Lewis for bringing our attention to the semialgebraic maps and for his decisive contribution in proving Theorem~\ref{threeormoreplayers}, and Prof. S. Sorin for bringing our attention to the paper \cite{doi:10.1287/moor.11.1.1}.

\biboptions{sort}
\bibliographystyle{plain}
\bibliography{ref}

\appendix

\section{Proof of Proposition~\ref{prop:dominatedoutcomes3players}}\label{appendix}

We need to prove that neither the  coalitions made by two players nor the grand coalition have incentive to deviate from the Nash equilibrium. First of, the coalitions made by two players.
We fix the strategy of player 3, but since the game is completely
symmetric, the same argument holds for every player. We calculate
the expected utility of player~1 and player~2 and we force them to be strictly greater
than $\frac{1}{2}$.

\begin{eqnarray*}
2\cdot\frac{1}{2}x_{11}x_{21}+2\cdot\frac{1}{2}x_{12}x_{22} & > & \frac{1}{2}\\
2\cdot\frac{1}{2}x_{11}x_{22}+2\cdot\frac{1}{2}x_{12}x_{21} & > & \frac{1}{2}
\end{eqnarray*}
Since $x_{i1}+x_{i2}=1$ we can write:
\begin{eqnarray*}
x_{11}x_{21}+\left(1-x_{11}\right)\left(1-x_{21}\right) & > & \frac{1}{2}\\
x_{11}\left(1-x_{21}\right)+\left(1-x_{11}\right)x_{21} & > & \frac{1}{2}
\end{eqnarray*}
and we obtain:
\begin{eqnarray*}
2x_{11}x_{21}-x_{11}-x_{21}+1 & > & \frac{1}{2}\\
x_{11}+x_{21}-2x_{11}x_{21} & > & \frac{1}{2}
\end{eqnarray*}

By summing the two constraints we obtain $1>1$ and therefore the
system is unfeasible.
We do now the calculations for the grand coalition.
\begin{eqnarray*}
2x_{11}x_{21}x_{31}+2x_{12}x_{22}x_{32} & > & \frac{1}{2}\\
2x_{11}x_{22}x_{31}+2x_{12}x_{21}x_{32} & > & \frac{1}{2}\\
2x_{11}x_{22}x_{32}+2x_{12}x_{21}x_{31} & > & \frac{1}{2}
\end{eqnarray*}
By replacing $x_{i2}$ by $1-x_{i1}$ we can write:
\begin{eqnarray*}
2x_{11}x_{21}x_{31}+2\left(1-x_{11}\right)\left(1-x_{21}\right)\left(1-x_{31}\right) & > & \frac{1}{2}\\
2x_{11}\left(1-x_{21}\right)x_{31}+2\left(1-x_{11}\right)x_{21}\left(1-x_{31}\right) & > & \frac{1}{2}\\
2x_{11}\left(1-x_{21}\right)\left(1-x_{31}\right)+2\left(1-x_{11}\right)x_{21}x_{31} & > & \frac{1}{2}
\end{eqnarray*}
and we obtain:
\begin{eqnarray*}
2-2x_{21}-2x_{31}+2x_{21}x_{31}-2x_{11}+2x_{11}x_{21}+2x_{11}x_{31} & > & \frac{1}{2}\\
2x_{11}x_{31}+2x_{21}-2x_{11}x_{21}-2x_{21}x_{31} & > & \frac{1}{2}\\
2x_{11}-2x_{11}x_{21}-2x_{11}x_{31}+2x_{21}x_{31} & > & \frac{1}{2}
\end{eqnarray*}

\textcolor{black}{By isolating $x_{11}$, we can write the following system of inequalities:
\begin{align}
x_{11}(x_{31}+x_{21}-1)  &  >      \frac{1}{4}-(1-x_{21})(1-x_{31})  \label{masterinequality1}    \\
x_{11}(x_{31}-x_{21})    &  >      \frac{1}{4}-x_{21}(1-x_{31})      \label{masterinequality2}    \\
x_{11}(1-x_{31}-x_{21})  &  >      \frac{1}{4}-x_{21} x_{31}         \label{masterinequality3}    \\
x_{11},x_{21},x_{31}     &  \geq   0                                   \nonumber  \\
x_{11},x_{21},x_{31}     &  \leq   1                                   \nonumber
\end{align}}

\textcolor{black}{We show that this inequality system is never satisfied. To do this, we consider all the possible cases characterized by the sign ($=,<,>0$) of $x_{31}+x_{21}-1$ and $x_{31}-x_{21}$.}

\textcolor{black}{The case in which $x_{31}+x_{21}-1 = 0$, independently of the value of $x_{31}-x_{21}$, inequality~(\ref{masterinequality3}) can be written as $\frac{1}{4} - x_{21} (1-x_{21})>0$ that is never satisfied given that the maximum  of $\frac{1}{4} - x_{21} (1-x_{21})$ is $0$ for $x_{21}=\frac{1}{2}$. The case in which $x_{31}-x_{21} = 0$, independently of the value of $1-x_{31}-x_{21}$, inequality~(\ref{masterinequality2}) can be written as $\frac{1}{4} - x_{21} (1-x_{21})>0$ that, as proved before, is never satisfied.}

\textcolor{black}{The case in which $(x_{31}+x_{21}-1 > 0~\wedge~x_{31}-x_{21}>0)$ we can write inequalities~(\ref{masterinequality1}), (\ref{masterinequality2}), (\ref{masterinequality3}) as
\begin{align}
x_{11}    &  >      \frac{\frac{1}{4}-(1-x_{21})(1-x_{31})}{x_{31}+x_{21}-1}  \label{inequality1}    \\
x_{11}    &  >      \frac{\frac{1}{4}-x_{21}(1-x_{31})}{x_{31}-x_{21}}        \label{inequality2}    \\
x_{11}    &  <      \frac{x_{2,1} x_{31}-\frac{1}{4}}{x_{31}+x_{21}-1}         \label{inequality3}
\end{align}}

\textcolor{black}{By combining inequalities~(\ref{inequality1}) and~(\ref{inequality3}) we obtain:
\begin{align}
\frac{\frac{1}{4}-(1-x_{21})(1-x_{31})}{x_{31}+x_{21}-1} - \frac{x_{21} x_{31}-\frac{1}{4}}{x_{31}+x_{21}-1} &<0\nonumber \\
x_{21}+x_{31}-2x_{21}x_{31} -\frac{1}{2}&<0      \label{inequality4}
\end{align}}

\textcolor{black}{By combining inequalities~(\ref{inequality2}) and~(\ref{inequality3}) we obtain:
\begin{align}
\frac{\frac{1}{4}-x_{21}(1-x_{31})}{x_{31}-x_{21}} - \frac{x_{21} x_{31}-\frac{1}{4}}{x_{31}+x_{21}-1} &<0\nonumber \\
x_{21}   & \neq \frac{1}{2}   \nonumber\\
x_{31}   & <  \frac{1}{2}  \label{inequality5}
\end{align}}

\textcolor{black}{From $x_{31}-x_{21}>0$ and inequality~(\ref{inequality5}), it follows that $x_{21}x_{31}<\frac{1}{4}$. From this last inequality and inequality~(\ref{inequality4}), it follows that $x_{21}+x_{31}-1<0$. Given that we are assuming $x_{31}+x_{21}-1 > 0$, the inequality system cannot be satisfied.}

\textcolor{black}{The cases in which $(x_{31}+x_{21}-1 > 0 ~\wedge ~x_{31}-x_{21}<0)$, $(x_{31}+x_{21}-1 < 0~\wedge~x_{31}-x_{21}>0)$, and $(x_{31}+x_{21}-1 > 0~\wedge~x_{31}-x_{21}<0)$ lead to calculations similar to the case in which $(x_{31}+x_{21}-1 > 0~\wedge~x_{31}-x_{21}>0)$. This completes the proof.}\newendproof

\end{document}